# Continuous dependence on data for linear evolution PDEs on the quarter-plane

Andreas Chatziafratis [1,3,*] and Spyridon Kamvissis [2,3]

**Abstract.** In this note, we announce a systematic analysis of "continuous dependence on the data" in classical spaces for the initial-boundary-value problem of the diffusion equation on the half-line, with data that are not necessarily compatible at the quadrant corner. This is based on a recent approach [1-9] to rigorously analyzing integral representations derived via the unified transform method [10-18] of Fokas. No exotic phenomena were discovered in this case, yet our findings appear to be new in the pertinent literature. These results supplement our previous investigations [2,6] on existence and (non)uniqueness within the framework of "well-posedness". The present detailed exposition elucidates the subtleties involved while also demonstrating a generic technique. Applications of the latter to several other IBVPs and PDEs will be reported elsewhere.

***Convergence in the spaces*** $\mathcal{S}([0,\infty))$, $\mathcal{S}_1(\overline{Q} - \{(0,0)\})$, $C^\infty([0,\infty))$, $\mathcal{S}_1(\overline{Q})$.

**1.** Let $\mathcal{S}([0,\infty))$ be the space of *Schwartz* functions in the half-line. More precisely,

$$\mathcal{S}([0,\infty)) := \{u \in C^\infty([0,\infty)) : \rho_N(u) < \infty \text{ for every nonnegative integer } N\}$$

where

$$\rho_N(u) := \sup\left\{(1+x)^N \left|\frac{\partial^k u(x)}{\partial x^k}\right| : k \leq N, x \geq 0\right\}.$$

Endowed with the topology defined by the norms $\{\rho_N : N \in \mathbb{N} \cup \{0\}\}$, $\mathcal{S}([0,\infty))$ becomes a Fréchet space. We recall that this means that the topology of $\mathcal{S}([0,\infty))$ is defined by the metric

$$\rho(u,v) := \sum_N \frac{1}{2^N} \frac{\rho_N(u-v)}{1+\rho_N(u-v)}, \quad u,v \in \mathcal{S}([0,\infty)).$$

Thus, a sequence $u_s \to u$, in $\mathcal{S}([0,\infty))$ with respect to the metric $\rho$, as $s \to \infty$, if and only if $\rho_N(u_s - u) \to 0$, as $s \to \infty$, for every $N$. Equivalently, $u_s \to u$ if and only if

$$x^L \frac{\partial^k u_s(x)}{\partial x^k} \to x^L \frac{\partial^k u(x)}{\partial x^k}, \text{ uniformly for } x \geq 0 \text{ and for all nonnegative integers } k \text{ and } L.$$

**2.** Let $C^\infty(\overline{Q} - \{(0,0)\})$ be the set of the functions $h = h(x,t)$ which are $C^\infty$ for $(x,t) \in Q - \{(0,0)\}$, such that the limits

$$\left.\frac{\partial^{k+l} h(x,t)}{\partial x^k \partial t^l}\right|_{(x,t)=(0,\tau)} := \lim_{Q \ni (x,t) \to (0,\tau)} \frac{\partial^{k+l} h(x,t)}{\partial x^k \partial t^l}, \text{ for } \tau > 0,$$

and

$$\left.\frac{\partial^{k+l} h(x,t)}{\partial x^k \partial t^l}\right|_{(x,t)=(\chi,0)} := \lim_{Q \ni (x,t) \to (\chi,0)} \frac{\partial^{k+l} h(x,t)}{\partial x^k \partial t^l}, \text{ for } \chi > 0,$$

exist for all nonnegative integers $k,l$, and the functions

$$\left.\frac{\partial^{k+l} h(x,t)}{\partial x^k \partial t^l}\right|_{(x,t)=(0,\tau)} (\tau > 0), \quad \left.\frac{\partial^{k+l} h(x,t)}{\partial x^k \partial t^l}\right|_{(x,t)=(\chi,0)} (\chi > 0)$$

are $C^\infty$ with respect to $\tau > 0$ and $\chi > 0$, respectively.

For example, if $h \in C^\infty(O)$, where $O \subset \mathbb{R}^2$ is an open set with $O \supset \overline{Q} - \{(0,0)\}$, then $h \in C^\infty(\overline{Q} - \{(0,0)\})$.

---

[1] Department of Mathematics, National and Kapodistrian University of Athens

[2] Department of Pure and Applied Mathematics, University of Crete

[3] Institute of Applied and Computational Mathematics, FORTH, Crete

*corresponding author, e-mail: chatziafrati@math.uoa.gr

**3.** For $n \in \mathbb{N}$, we consider the sets

$$\Lambda_n := \{(x,t) \in \overline{Q} : \tfrac{1}{n} \leq t \leq n\} \cup \{(x,t) \in \overline{Q} : x \geq \tfrac{1}{n}, t \leq \tfrac{1}{n}\},$$

which exhaust $\overline{Q} - \{(0,0)\}$, i.e., $\overline{Q} - \{(0,0)\} = \bigcup_{n \in \mathbb{N}} \Lambda_n$. (See Fig 1)

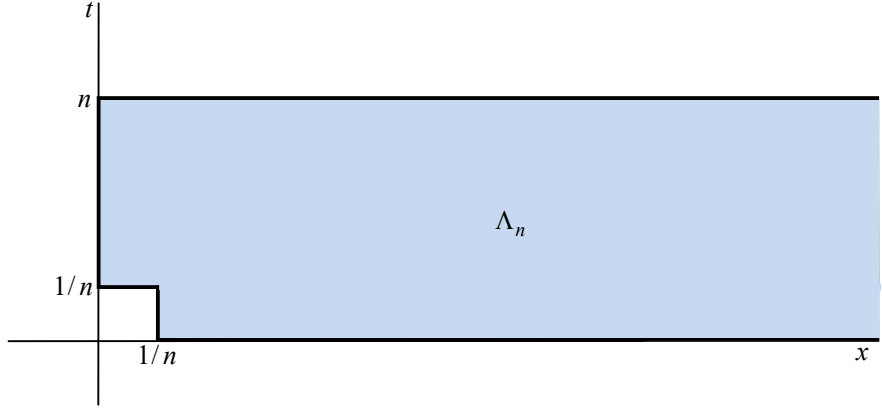

**Fig 1** The infinite closed "strips" $\Lambda_n = \{x \geq 0, \tfrac{1}{n} \leq t \leq n\} \cup \{x \geq \tfrac{1}{n}, 0 \leq t \leq \tfrac{1}{n}\}$.

Then we define

$$\lambda_n(h) := \sup\left\{(1+x)^n \left|\frac{\partial^{k+l} h(x,t)}{\partial x^k \partial t^l}\right| : (x,t) \in \Lambda_n, k+l \leq n\right\}, \text{ for } h \in C^\infty(\overline{Q} - \{(0,0)\}) \text{ and } n \in \mathbb{N},$$

and we let

$$\mathcal{S}_1(\overline{Q} - \{(0,0)\}) = \{h \in C^\infty(\overline{Q} - \{(0,0)\}) : \lambda_n(h) < \infty \text{ for every } n \in \mathbb{N}\}.$$

We endow $\mathcal{S}_1(\overline{Q} - \{(0,0)\})$ with the topology induced by the seminorms $\{\lambda_n : n \in \mathbb{N}\}$. Thus $\mathcal{S}_1(\overline{Q} - \{(0,0)\})$ becomes a metric space with the metric

$$\lambda(h_1, h_2) := \sum_{n \in \mathbb{N}} \frac{1}{2^n} \frac{\lambda_n(h_1 - h_2)}{1 + \lambda_n(h_1 - h_2)}, \quad h_1, h_2 \in \mathcal{S}_1(\overline{Q} - \{(0,0)\}).$$

**4.** Similarly, exhausting $\overline{Q}$ with the sets (Fig 2)

$$\Lambda_n^* := \{(x,t) \in \overline{Q} : t \leq n\}, \quad n \in \mathbb{N},$$

we define

$$\mathcal{S}_1(\overline{Q}) = \{f \in C^\infty(\overline{Q}) : \lambda_n^*(f) < \infty \text{ for every } n \in \mathbb{N}\},$$

where

$$\lambda_n^*(f) := \sup\left\{(1+x)^n \left|\frac{\partial^{k+l} f(x,t)}{\partial x^k \partial t^l}\right| : (x,t) \in \Lambda_n^*, k+l \leq n\right\},$$

and we endow $\mathcal{S}_1(\overline{Q})$ with the topology induced by the seminorms $\{\lambda_n^* : n \in \mathbb{N}\}$.

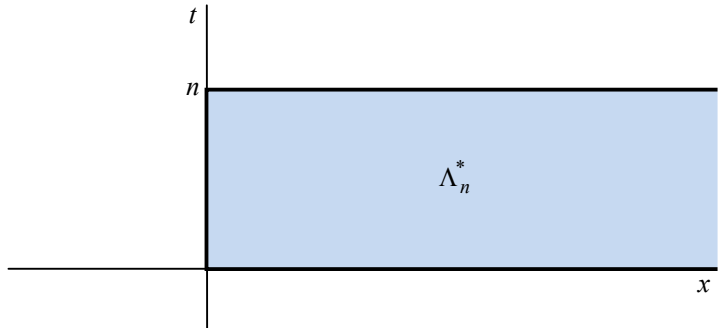

**Fig. 2** The strips $\Lambda_n^*$

**5.** The topology in $C^\infty([0,\infty))$ is defined by the seminorms

$$\sup\left\{\left|\frac{\partial^l g(t)}{\partial x^l}\right|: 0\le t\le n\right\}, \quad g\in C^\infty([0,\infty)),\ l\in\mathbb{N}\cup\{0\},\ n\in\mathbb{N}.$$

A sequence $g_s \to g$, in $C^\infty([0,\infty))$, as $s\to\infty$, if and only if

$$\frac{\partial^l g_s(t)}{\partial t^l}\to\frac{\partial^l g(t)}{\partial t^l},\text{ uniformly for } 0\le t\le n,\text{ for all nonnegative integers } l \text{ and for all } n\in\mathbb{N}.$$

**Lemma 1** *Let* $\mathbb{I}_o^+ : \mathcal{S}([0,\infty))\to \mathcal{S}_1(\overline{Q}-\{(0,0)\})$ *be defined as follows: For* $u\in\mathcal{S}([0,\infty))$,

$$(\mathbb{I}_o^+ u)(x,t)=\int_{-\infty}^\infty e^{i\lambda x-\lambda^2 t}\hat{u}(\lambda)d\lambda,\text{ for }(x,t)\in Q,$$

$$(\mathbb{I}_o^+ u)(0,t)=\int_{-\infty}^\infty e^{-\lambda^2 t}\hat{u}(\lambda)d\lambda,\text{ for }t>0,$$

$$(\mathbb{I}_o^+ u)(x,0)=\int_{-1}^1 e^{i\lambda x}\hat{u}(\lambda)d\lambda+u(0)\int_{\gamma^*}e^{i\lambda x}\frac{d\lambda}{i\lambda}+\left(\int_{-\infty}^{-1}+\int_1^\infty\right)e^{i\lambda x}\frac{\widehat{(u')}(\lambda)}{i\lambda}d\lambda,\text{ for }x>0,$$

*where* $\gamma^*:=(\gamma\cap\{|\lambda|\ge\sqrt{2}\})+[-1+i,-1]+[1,1+i]$ (Fig 3).
*Then* $\mathbb{I}_o^+$ *is well-defined and maps* $\mathcal{S}([0,\infty))$ *continuously to* $\mathcal{S}_1(\overline{Q}-\{(0,0)\})$.

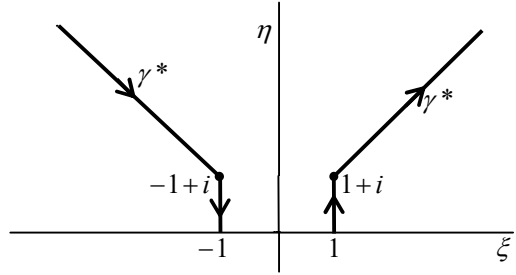

**Fig. 3** *The contour* $\gamma^*$

**Proof** First, the integrals, which define $(\mathbb{I}_o^+ u)(x,t)$ for the various $(x,t)\in\overline{Q}-\{(0,0)\}$, are absolutely convergent.

Let $u_s\in\mathcal{S}([0,\infty))$, $s=1,2,3,...$, be a sequence so that $u_s\to 0$, in $\mathcal{S}([0,\infty))$, as $s\to\infty$. Let us also fix $k,l$.

***Step 1*** Setting $\Lambda_n^{(1)}=\{x\ge 0,\frac{1}{n}\le t\le n\}$ (fig. 4) and $\Lambda_n^{(2)}=\{x\ge\frac{1}{n},0\le t\le n\}$ (fig. 5), we have

$$(\mathbb{I}_o^+ u)(x,t)=\int_{-\infty}^\infty e^{i\lambda x-\lambda^2 t}\hat{u}(\lambda)d\lambda,\text{ for }(x,t)\in\Lambda_n^{(1)},$$

and

$$(\mathbb{I}_o^+ u)(x,t)=\int_{-1}^1 e^{i\lambda x-\lambda^2 t}\hat{u}(\lambda)d\lambda+\sum_{j=1}^M u^{(j-1)}(0)\int_{\gamma^*}e^{i\lambda x-\lambda^2 t}\frac{d\lambda}{(i\lambda)^j}+\left(\int_{-\infty}^{-1}+\int_1^\infty\right)\left\{e^{i\lambda x-\lambda^2 t}[\widehat{u^{(M)}}](\lambda)\frac{d\lambda}{(i\lambda)^M}\right\},\text{ for }(x,t)\in\Lambda_n^{(2)}.$$

Differentiating we obtain

$$\frac{\partial^{k+l}(\mathbb{I}_o^+ u)(x,t)}{\partial x^k\partial t^l}=\int_{-\infty}^\infty e^{i\lambda x-\lambda^2 t}(i\lambda)^k(-\lambda^2)^l\hat{u}(\lambda)d\lambda,\text{ for }(x,t)\in\Lambda_n^{(1)},\qquad(1)$$

and, for $M>k+2l+1$,

$$\frac{\partial^{k+l}(\mathbb{I}_o^+ u)(x,t)}{\partial x^k\partial t^l}=\int_{-1}^1 e^{i\lambda x-\lambda^2 t}(i\lambda)^k(-\lambda^2)^l\hat{u}(\lambda)d\lambda+\sum_{j=1}^M u^{(j-1)}(0)\int_{\gamma^*}e^{i\lambda x-\lambda^2 t}(i\lambda)^k(-\lambda^2)^l\frac{d\lambda}{(i\lambda)^j}$$

$$+\left(\int_{-\infty}^{-1}+\int_{1}^{\infty}\right)\left\{e^{i\lambda x-\lambda^{2}t}(i\lambda)^{k}(-\lambda^{2})^{l}[u^{(M)}\hat{\,}](\lambda)\frac{d\lambda}{(i\lambda)^{M}}\right\}, \text{ for } (x,t)\in\Lambda_{n}^{(2)}, \quad (2)$$

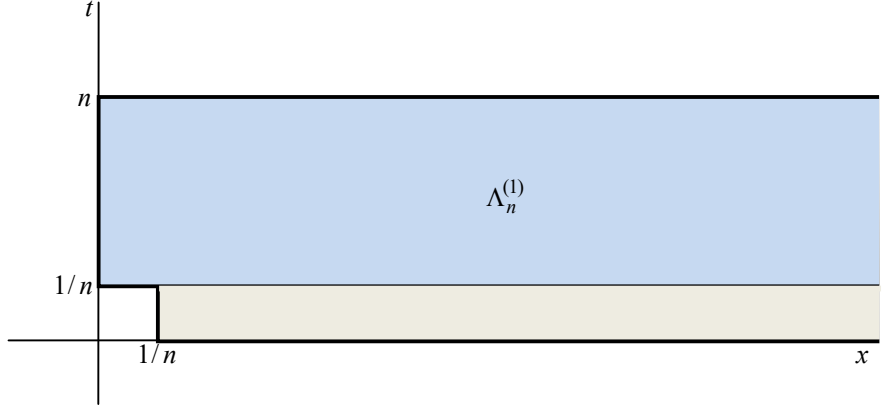

**Fig 4** The strips $\Lambda_{n}^{(1)} = \{x \geq 0, \frac{1}{n} \leq t \leq n\}$.

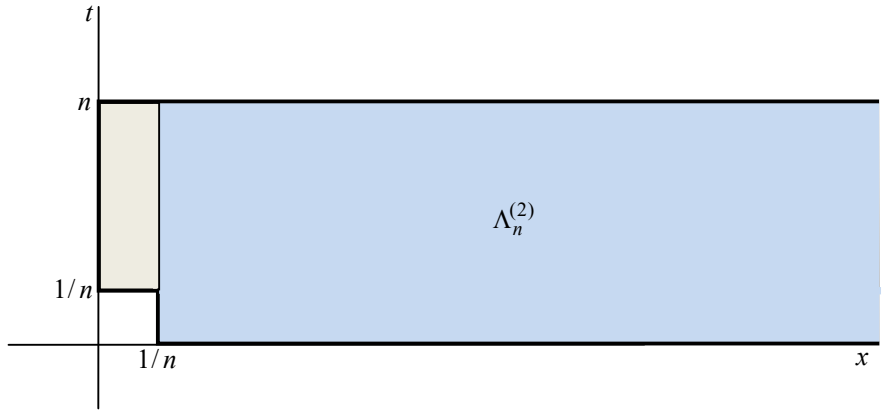

**Fig 5** The strips $\Lambda_{n}^{(2)} = \{x \geq \frac{1}{n}, 0 \leq t \leq n\}$.

It follows that $\mathbb{I}_o^+ u \in C^{\infty}(\overline{Q} - \{(0,0)\})$, and, in view of Theorem 3, we obtain that, furthermore, $\mathbb{I}_o^+ u \in \mathcal{S}_1(\overline{Q} - \{(0,0)\})$, i.e., $\mathbb{I}_o^+ : \mathcal{S}([0,\infty)) \to \mathcal{S}_1(\overline{Q} - \{(0,0)\})$ is well-defined.

**Step 2** Let us fix $M$, with $M > k + 2l + 1$. Then for every $\varepsilon > 0$, there is $s_0 = s_0(\varepsilon)$ so that

$$s \geq s_0 \implies (1+y)^2 \left| \frac{d^M u_s(y)}{dy^M} \right| \leq \varepsilon, \text{ for every } y \geq 0.$$

It follows that

$$s \geq s_0 \implies |[u_s^{(M)}\hat{\,}](\lambda)| = \left|\int_0^{\infty} e^{-i\lambda y} \frac{d^M u_s(y)}{dy^M} dy\right| \leq \varepsilon \int_0^{\infty} \frac{dy}{(1+y)^2} = \varepsilon, \text{ for every } \lambda \in \mathbb{R}.$$

Therefore,

$$\sup_{(x,t)\in\overline{Q}}\left|\left(\int_{-\infty}^{-1}+\int_{1}^{\infty}\right)\left\{e^{i\lambda x-\lambda^{2}t}(i\lambda)^{k}(-\lambda^{2})^{l}[u_s^{(M)}\hat{\,}](\lambda)\frac{d\lambda}{(i\lambda)^{M}}\right\}\right| \leq \varepsilon\left(\int_{-\infty}^{-1}+\int_{1}^{\infty}\right)\frac{d\lambda}{|\lambda|^{M-k-2l}} \leq 2\varepsilon, \text{ for } s \geq s_0. \quad (3)$$

Thus, the LHS of (3) converges to $0$, as $s \to \infty$.
Similarly,

$$\sup_{(x,t)\in\overline{Q}}\left|\int_{-1}^{1} e^{i\lambda x-\lambda^{2}t}(i\lambda)^{k}(-\lambda^{2})^{l}\hat{u}_s(\lambda)d\lambda\right| \to 0, \text{ as } s \to \infty, \quad (4)$$

and, for fixed $n$,

$$\sup_{(x,t)\in\Lambda_n^{(1)}}\left|\int_{-\infty}^{\infty}e^{i\lambda x-\lambda^2 t}(i\lambda)^k(-\lambda^2)^l\hat{u}_s(\lambda)d\lambda\right|\to 0, \text{ as } s\to\infty. \tag{5}$$

**Step 3** With $u_s$ as in *Step 2*, $u_s^{(j-1)}(0)\to 0$, for $1\le j\le M$, and therefore, for fixed $n$,

$$\sup_{(x,t)\in\Lambda_n^{(2)}}\left|\sum_{j=1}^{M}u_s^{(j-1)}(0)\int_{\gamma^*}e^{i\lambda x-\lambda^2 t}(i\lambda)^k(-\lambda^2)^l\frac{d\lambda}{(i\lambda)^j}\right|\to 0, \text{ as } s\to\infty. \tag{6}$$

**Step 4** In view of (2), (3), (4) and (6),

$$\sup_{(x,t)\in\Lambda_n^{(2)}}\left|\frac{\partial^{k+l}(\mathbb{I}_\circ^+ u)(x,t)}{\partial x^k \partial t^l}\right|\to 0, \text{ as } s\to\infty. \tag{7}$$

**Step 5** Multipling (1) by $(ix)^L$ and integrating by parts, we obtain

$$(ix)^L\frac{\partial^{k+l}(\mathbb{I}_\circ^+ u)(x,t)}{\partial x^k \partial t^l}=(-1)^L\int_{-\infty}^{\infty}e^{i\lambda x}\frac{d^L}{d\lambda^L}[e^{-\lambda^2 t}(i\lambda)^k(-\lambda^2)^l\hat{u}(\lambda)]d\lambda, \text{ for } (x,t)\in\Lambda_n^{(1)}. \tag{8}$$

It follows that, for fixed $L$, $k$, $l$ and $n$,

$$\sup_{(x,t)\in\Lambda_n^{(1)}}\left|x^L\frac{\partial^{k+l}(\mathbb{I}_\circ^+ u_s)(x,t)}{\partial x^k \partial t^l}\right|\to 0, \text{ as } s\to\infty. \tag{9}$$

Similarly, writing (2) with *sufficiently* large $M$, multiplying it by $(ix)^L$, substituting $(ix)^L e^{i\lambda x}$ by $d^L(e^{i\lambda x})/d\lambda^L$ in the integrals of the RHS of the resulting (2), and integrating by parts, we see that, for fixed $L$, $k$, $l$ and $n$,

$$\sup_{(x,t)\in\Lambda_n^{(2)}}\left|x^L\frac{\partial^{k+l}(\mathbb{I}_\circ^+ u_s)(x,t)}{\partial x^k \partial t^l}\right|\to 0, \text{ as } s\to\infty. \tag{10}$$

Indeed, the proof of (9) is analogous to the proof (7).

*Completion of the proof* The conclusion of the lemma follows from the conclusions of *Steps 1-5*.

**Lemma 2** *Let* $\mathbb{I}_1:C^\infty([0,\infty))\to\mathcal{S}_1(\overline{Q}-\{(0,0)\})$ *be defined as follows: For* $g\in C^\infty([0,\infty))$,

$$(\mathbb{I}_1 g)(x,t)=\int_\gamma e^{i\lambda x-\lambda^2 t}\tilde{g}(\omega(\lambda),t)\lambda d\lambda, \text{ for } (x,t)\in\overline{Q}-\{(0,0)\}. \text{ (In this section, } \omega(\lambda)=\lambda^2.)$$

*Then* $\mathbb{I}_1$ *is well-defined and maps* $C^\infty([0,\infty))$ *continuously to* $\mathcal{S}_1(\overline{Q}-\{(0,0)\})$.

**Proof** First, the integral which defines $(\mathbb{I}_1 g)(x,t)$ converges absolutely, when $x>0$, and it is equal to

$$\int_{-\infty}^{\infty}e^{i\lambda x-\lambda^2 t}\tilde{g}(\omega(\lambda),t)\lambda d\lambda$$

when $(x,t)\in Q$.

**Step 1** We have

$$(\mathbb{I}_1 g)(x,t)=\int_{\gamma_0}e^{i\lambda x-\lambda^2 t}\tilde{g}(\omega(\lambda),t)\lambda d\lambda.$$

(Recall $\gamma_0=(\gamma\cap\{|\lambda|\ge\sqrt{2}\})+[-1+i,1+i]$, see also Fig 2 in [6])
First, let us keep in mind that

$$\left|e^{i\lambda x}\right|=e^{-[x\sin(\pi/4)]|\lambda|} \text{ and } x^L|\lambda|^N e^{i\lambda x}\le x^L|\lambda|^N e^{-[x\sin(\pi/4)]|\lambda|}, \text{ for } \lambda\in\gamma. \tag{11}$$

In order to deal with the integral over $\gamma\cap\{|\lambda|\ge\sqrt{2}\}$, we use the inequality

$$\left|e^{-\omega(\lambda)t}\tilde{g}(\omega(\lambda),t)\right|\le\int_0^t|g(\tau)|d\tau, \text{ for } \lambda\in\gamma,$$

in combination with (11), to obtain that, for $x \geq 1/n$,

$$x^L \left| \int_{\gamma \cap \{|\lambda| \geq 1\}} \lambda^N e^{i\lambda x - \lambda^2 t} \lambda \widetilde{g}(\omega(\lambda),t) d\lambda \right| \leq \left( \int_0^t |g(\tau)| d\tau \right) x^L e^{-[x\sin(\pi/4)]/2} \int_{\gamma \cap \{|\lambda| \geq 1\}} |\lambda|^{N+1} e^{-[\sin(\pi/4)]|\lambda|/(2n)} d|\lambda|. \quad (12)$$

On the other hand, for $\lambda \in [-1+i, 1+i]$ we have $|\lambda| \leq \sqrt{2}$, $\mathrm{Im}\,\lambda = 1$ and $|e^{i\lambda x}| = e^{-x\,\mathrm{Im}\,\lambda} = e^{-x}$, and, therefore,

$$x^L \left| \int_{[-1+i,1+i]} \lambda^L e^{i\lambda x - \lambda^2 t} \lambda \widetilde{g}(\omega(\lambda),t) d\lambda \right| \leq \left( \int_0^t |g(\tau)| d\tau \right) x^L e^{-x} \int_{[-1+i,1+i]} d|\lambda|.$$

Now,

$$x^L \frac{\partial^k (\mathbb{I}_1 g)(x,t)}{\partial x^k} = x^L \int_{\gamma^*} e^{i\lambda x - \lambda^2 t} (i\lambda)^k \widetilde{g}(\omega(\lambda),t) \lambda d\lambda,$$

in combination with (12), implies that, for fixed $L, k$ and $n$,

$$\sup_{(x,t) \in \Lambda_n^{(2)}} \left| x^L \frac{\partial^k (\mathbb{I}_1 g)(x,t)}{\partial x^k} \right| \preceq \sup_{|t| \leq n} |g(t)|.$$

Similar computations show that, more generally, for fixed $L, k, l$ and $n$,

$$\sup_{(x,t) \in \Lambda_n^{(2)}} \left| x^L \frac{\partial^{k+l} (\mathbb{I}_1 g)(x,t)}{\partial x^k \partial t^l} \right| \preceq \sup\{ |g^{(\ell)}(t)| : |t| \leq n, \ell \leq l \} \quad (\forall g).$$

**Step 2** Let $g(t) \in C^\infty([0,\infty))$. For $\lambda \in \mathbb{C}, \lambda \neq 0$, we have

$$e^{-\lambda^2 t} \int_{\tau=0}^t e^{\lambda^2 \tau} g(\tau) d\tau = \frac{g(t)}{\lambda^2} - \frac{g(0)}{\lambda^2} e^{-\lambda^2 t} - \frac{1}{\lambda^2} e^{-\lambda^2 t} \int_{\tau=0}^t e^{\lambda^2 \tau} g'(\tau) d\tau$$

$$= \frac{g(t)}{\lambda^2} - \frac{g(0)}{\lambda^2} e^{-\lambda^2 t} - \frac{g'(t)}{\lambda^4} + \frac{g'(0)}{\lambda^4} e^{-\lambda^2 t} + \frac{1}{\lambda^4} e^{-\lambda^2 t} \int_{\tau=0}^t e^{\lambda^2 \tau} g''(\tau) d\tau$$

$$= \sum_{j=0}^M (-1)^n \left[ \frac{g^{(j)}(t)}{\lambda^{2j+2}} - \frac{g^{(j)}(0)}{\lambda^{2j+2}} e^{-\lambda^2 t} \right] + (-1)^{M+1} \frac{1}{\lambda^{2M+2}} e^{-\lambda^2 t} \int_{\tau=0}^t e^{\lambda^2 \tau} g^{(M+1)}(\tau) d\tau. \quad (13)$$

Setting $\gamma_1 := \gamma \cap (\{|\lambda| \leq \sqrt{2}\})$ and $\gamma_2 := \gamma \cap (\{|\lambda| \geq \sqrt{2}\})$ (Figs 6 & 7) we have

$$(\mathbb{I}_1 g)(x,t) = \int_{\gamma_1} e^{i\lambda x - \lambda^2 t} \widetilde{g}(\omega(\lambda),t) \lambda d\lambda + \int_{\gamma_2} e^{i\lambda x - \lambda^2 t} \widetilde{g}(\omega(\lambda),t) \lambda d\lambda$$

$$= \int_{\gamma_1} e^{i\lambda x - \lambda^2 t} \widetilde{g}(\omega(\lambda),t) \lambda d\lambda + \sum_{j=0}^M (-1)^j \left[ g^{(j)}(t) \int_{\gamma_2} e^{i\lambda x} \frac{d\lambda}{\lambda^{2j+1}} - g^{(j)}(0) \int_{\gamma_2} e^{i\lambda x - \lambda^2 t} \frac{d\lambda}{\lambda^{2j+1}} \right]$$

$$+ (-1)^{M+1} \int_{\gamma_2} e^{i\lambda x - \lambda^2 t} \left[ \int_{\tau=0}^t e^{\lambda^2 \tau} g^{(M+1)}(\tau) d\tau \right] \frac{d\lambda}{\lambda^{2M+1}}$$

$$= \int_{\gamma_1} e^{i\lambda x - \lambda^2 t} \widetilde{g}(\omega(\lambda),t) \lambda d\lambda$$

$$+ \sum_{j=0}^M (-1)^j \left[ -g^{(j)}(t) \int_{[-1+i,\,1+i]} e^{i\lambda x} \frac{d\lambda}{\lambda^{2j+1}} - g^{(j)}(0) \left( \int_{-\infty}^{-1} + \int_1^\infty \right) e^{i\lambda x - \lambda^2 t} \frac{d\lambda}{\lambda^{2j+1}} - g^{(j)}(0) \left( \int_{[-1,-1+i]} + \int_{[1+i,1]} \right) e^{i\lambda x - \lambda^2 t} \frac{d\lambda}{\lambda^{2j+1}} \right]$$

$$+ (-1)^{M+1} \int_{\gamma_2} e^{i\lambda x - \lambda^2 t} \left[ \int_{\tau=0}^t e^{\lambda^2 \tau} g^{(M+1)}(\tau) d\tau \right] \frac{d\lambda}{\lambda^{2M+1}}. \quad (14)$$

**Step 3** Differentiating (14) and multiplying by $(ix)^\ell$, we obtain:

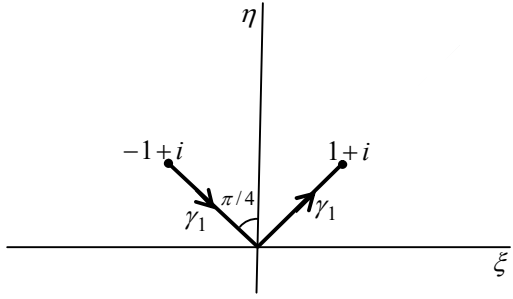
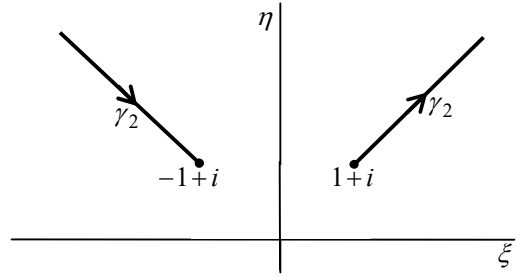

**Fig. 6** *The contour* $\gamma_1$  **Fig. 7** *The contour* $\gamma_2$

$$(ix)^L \frac{\partial^{k+l}(\mathbb{I}_1 g)(x,t)}{\partial x^k \partial t^l} = \int_{\gamma_1} \frac{d^L(e^{i\lambda x})}{d\lambda^L}\left\{(i\lambda)^k \frac{\partial^l}{\partial t^l}\left[e^{-\lambda^2 t}\int_{\tau=0}^{t} e^{\lambda^2 \tau} g(\tau)d\tau\right]\lambda\right\}d\lambda$$

$$-\sum_{j=0}^{M}(-1)^j g^{(j+l)}(t) \int_{[-1+i,1+i]} \frac{d^L(e^{i\lambda x})}{d\lambda^L}\left[(i\lambda)^k \frac{1}{\lambda^{2j+1}}\right]d\lambda$$

$$-\sum_{j=0}^{M}(-1)^j g^{(j)}(0)\left(\int_{-\infty}^{-1}+\int_{1}^{\infty}\right) \frac{d^L(e^{i\lambda x})}{d\lambda^L}\left[e^{-\lambda^2 t}(i\lambda)^k (-\lambda^2)^l \frac{1}{\lambda^{2j+1}}\right]d\lambda$$

$$-\sum_{j=0}^{M}(-1)^j g^{(j)}(0)\left(\int_{[-1,-1+i]}+\int_{[1+i,1]}\right) \frac{d^L(e^{i\lambda x})}{d\lambda^L}\left[e^{-\lambda^2 t}(i\lambda)^k (-\lambda^2)^l \frac{1}{\lambda^{2j+1}}\right]d\lambda$$

$$+(-1)^{M+1}\int_{\gamma_2} \frac{d^L(e^{i\lambda x})}{d\lambda^L}\left\{(i\lambda)^k \frac{\partial^l}{\partial t^l}\left[e^{-\lambda^2 t}\int_{\tau=0}^{t} e^{\lambda^2 \tau} g^{(M+1)}(\tau)d\tau\right]\frac{1}{\lambda^{2M+1}}\right\}d\lambda. \quad (15)$$

***Step 4*** Integrating by parts in (15), we obtain that, for fixed $L, k, l$ and $n$,

$$\sup_{(x,t)\in\Lambda_n^{(1)}}\left|x^L \frac{\partial^{k+l}(\mathbb{I}_1 g)(x,t)}{\partial x^k \partial t^l}\right| \preceq \sup\{|g^{(\ell)}(t)| : |t| \le n, \ell \le L+l\} \quad (\forall g). \quad (16)$$

*Completion of the proof* The conclusion of the lemma follows from the conclusions of *Steps 1-4*.

**Lemma 3** *Let* $\mathbb{I}_0^- : \mathcal{S}([0,\infty)) \to \mathcal{S}_1(\overline{Q} - \{(0,0)\})$ *be defined as follows: For* $u \in \mathcal{S}([0,\infty))$,

$$(\mathbb{I}_0^- u)(x,t) = \int_{\gamma} e^{i\lambda x - \lambda^2 t}\hat{u}(-\lambda)d\lambda, \text{ for } x > 0 \text{ and } t \ge 0,$$

$$(\mathbb{I}_0^- u)(x,t) = \int_{-\infty}^{\infty} e^{i\lambda x - \lambda^2 t}\hat{u}(-\lambda)d\lambda, \text{ for } x \ge 0 \text{ and } t > 0.$$

*Then* $\mathbb{I}_0^-$ *is well-defined and maps* $\mathcal{S}([0,\infty))$ *continuously to* $\mathcal{S}_1(\overline{Q} - \{(0,0)\})$.

**Proof** Let $u_s \to 0$, in $\mathcal{S}([0,\infty))$, as $s \to \infty$. Since for $(x,t) \in \Lambda_n^{(2)}$, $x \ge 1/n$, using the integral over $\gamma$, we prove, as in *Step 1* of Lemma 2 and *Step 2* of Lemma 1, that

$$\sup_{(x,t)\in\Lambda_n^{(2)}}\left|x^L \frac{\partial^{k+l}(\mathbb{I}_0^- u_s)(x,t)}{\partial x^k \partial t^l}\right| \to 0.$$

Also, since for $(x,t) \in \Lambda_n^{(1)}$, $t \ge 1/n$, using the integral over $(-\infty,\infty)$, we obtain

$$\sup_{(x,t)\in\Lambda_n^{(1)}}\left|x^L \frac{\partial^{k+l}(\mathbb{I}_0^- u_s)(x,t)}{\partial x^k \partial t^l}\right| \to 0.$$

The proofs of Lemmas 4 and 5, below, are similar to the above. To carry out these proofs we need also formulas (2.6), (2.8) and the analogue of (2.7) in [6].

**Lemma 4** *Let* $\mathbb{I}_2^+ : \mathcal{S}_1(\overline{Q}) \to \mathcal{S}_1(\overline{Q} - \{(0,0)\})$ *be defined as follows: For* $f \in \mathcal{S}_1(\overline{Q})$,

$$(\mathbb{I}_2^+ f)(x,t) = \int_{-\infty}^{\infty} e^{i\lambda x - \lambda^2 t} \widetilde{f}(\lambda, \omega(\lambda), t) d\lambda.$$

*Then* $\mathbb{I}_2^+$ *is well-defined and maps* $\mathcal{S}_1(\overline{Q})$ *continuously to* $\mathcal{S}_1(\overline{Q} - \{(0,0)\})$.

**Lemma 5** *Let* $\mathbb{I}_2^- : \mathcal{S}_1(\overline{Q}) \to \mathcal{S}_1(\overline{Q} - \{(0,0)\})$ *be defined as follows: For* $f \in \mathcal{S}_1(\overline{Q})$,

$$(\mathbb{I}_2^- f)(x,t) = \int_{\gamma} e^{i\lambda x - \lambda^2 t} \widetilde{f}(-\lambda, \omega(\lambda), t) d\lambda.$$

*Then* $\mathbb{I}_2^-$ *is well-defined and maps* $\mathcal{S}_1(\overline{Q})$ *continuously to* $\mathcal{S}_1(\overline{Q} - \{(0,0)\})$.

Finally, Lemmas 1-5 give the following theorem on continuous dependence on the data for Problem 2 of [6].

**Theorem** Let $\mathcal{D} := \mathcal{S}([0,\infty)) \times C^{\infty}([0,\infty)) \times \mathcal{S}_1(\overline{Q})$ be the space of the data $(u,g,f)$. Then the UTM solution operator $\Phi : \mathcal{D} \to \mathcal{S}_1(\overline{Q} - \{(0,0)\})$, of Problem 2, defined by

$$2\pi\Phi(u,g,f) = \mathbb{I}_0^+ u - \mathbb{I}_0^- u - 2i\mathbb{I}_1 g + \mathbb{I}_2^+ f - \mathbb{I}_2^- f, \text{ for } (u,g,f) \in \mathcal{D},$$

*is well-defined and continuous.*

***Comments 1.*** Let $\mathcal{D}_o := \{(u,g,f) \in \mathcal{D} : u(0) = g(0)\}$. Then $\mathcal{D}_o$ is a closed subspace of $\mathcal{D}$ and the restriction of the operator $\Phi$ to $\mathcal{D}_o$ maps $\mathcal{D}_o$ to $\mathcal{S}_1(\overline{Q} - \{(0,0)\}) \cap C^{(1,0)}([0,1] \times [0,1])$ and if a sequence $(u_s, g_s, f_s) \to (0,0,0)$ in $\mathcal{D}_o$, then $\Phi(u_s, g_s, f_s) \to 0$ both in $\mathcal{S}_1(\overline{Q} - \{(0,0)\})$ and $C^{(1,0)}([0,1] \times [0,1])$. (The space $C^{(1,0)}([0,1] \times [0,1])$ is the space of the continuous functions $h(x,t)$, $(x,t) \in [0,1] \times [0,1]$, which are also $C^1$ with respect to $x$, i.e., $\partial h(x,t) / \partial x$ is continuous for $(x,t) \in [0,1] \times [0,1]$. A sequence $h_s \to 0$ in $C^{(1,0)}([0,1] \times [0,1])$, if, both, $h_s(x,t) \to 0$ and $\partial h_s(x,t) / \partial x \to 0$, uniformly for $(x,t) \in [0,1] \times [0,1]$.)

**2.** Let $\mathcal{D}_1 := \{(u,g,f) \in \mathcal{D} : u(0) = g(0) \text{ and } g'(0) = u''(0) + f(0,0)\}$. Then $\mathcal{D}_1$ is a closed subspace of $\mathcal{D}$ and the restriction of the operator $\Phi$ to $\mathcal{D}_1$ maps $\mathcal{D}_1$ to $\mathcal{S}_1(\overline{Q} - \{(0,0)\}) \cap C^{(2,1)}([0,1] \times [0,1])$ and if the sequence $(u_s, g_s, f_s) \to (0,0,0)$ in $\mathcal{D}_1$, then $\Phi(u_s, g_s, f_s) \to 0$ both in $\mathcal{S}_1(\overline{Q} - \{(0,0)\})$ and $C^{(3,1)}([0,1] \times [0,1])$. (The space $C^{(3,1)}([0,1] \times [0,1])$ is the space of the $C^1$ functions $h(x,t)$, $(x,t) \in [0,1] \times [0,1]$, which are also $C^3$ with respect to $x$, i.e., $\partial^2 h(x,t) / \partial x^2$ and $\partial^3 h(x,t) / \partial x^3$ are continous for $(x,t) \in [0,1] \times [0,1]$. A sequence $h_s \to 0$ in $C^{(3,1)}([0,1] \times [0,1])$, if $h_s(x,t) \to 0$, $\partial h_s(x,t) / \partial x \to 0$, $\partial h_s(x,t) / \partial t \to 0$, $\partial^2 h_s(x,t) / \partial x^2 \to 0$ $\partial^3 h_s(x,t) / \partial x^3 \to 0$, and uniformly for $(x,t) \in [0,1] \times [0,1]$.)

**3.** Considering subspaces of $\mathcal{D}$ consisting of $(u,g,f)$ which satisfy higher order combatibility conditions at the origin, we obtain further continuity properties of the corresponding restrictions of the operator $\Phi$.

**4.** Results analogous to the above Theorem and remarks can be proved also for the UTM solution operator

$$(u,g,f) \to \mathbb{J}_0^+ u + \mathbb{J}_0^- u - \mathbb{J}_1 g + \mathbb{J}_2^+ f + \mathbb{J}_2^- f,$$

of *Problem 1* of [6]; see (1.5) therein.

**5.** More careful one should be with the corresponding problems for other equations such as the Airy equation on the negative half-line $U_t = U_{xxx} + f$ and the Schrödinger equation $iU_t = -U_{xx} + f$. For one thing, the solutions of the corresponding problems of these equations do not satisfy the analogues of Theorems 1 and 3 of [6], in general. The asymptotic behavior of these solutions, as $x \to \infty$, depends on certain compatibility conditions for the data at the origin. In addition, the limits of the solutions and their derivatives, as $t \to 0^+$, depend on similar compatibility conditions. (See [5,9].)

**6.** From the integral formula in the statement of the above Theorem, it is evident that the solution at $(x,t)$ depends only on past times $\tau \in [0,t)$ so one could readily state a Theorem with data $(g, f)$ as above but restricted to $[0,t)$.

# References

[1] A. Chatziafratis, Rigorous analysis of the Fokas method for linear evolution PDEs on the half-space, M.Sc. thesis, Advisors: N. Alikakos, G. Barbatis, I.G. Stratis, National and Kapodistrian University of Athens (2019). https://pergamos.lib.uoa.gr/uoa/dl/object/2877222
[2] A. Chatziafratis, Boundary behaviour of the solution of the heat equation on the half line via the unified transform method, *preprint* (2020). arXiv:2401.08331 [math.AP]
[3] A. Chatziafratis, D. Mantzavinos, Boundary behavior for the heat equation on the half-line, *Math. Methods Appl. Sci.* 45, 7364-93 (2022).
[4] A. Chatziafratis, S. Kamvissis, I.G. Stratis, Boundary behavior of the solution to the linear KdV equation on the half-line, *Stud. Appl. Math.* 150, 339-379 (2023).
[5] A. Chatziafratis, T. Ozawa, S.-F. Tian, Rigorous analysis of the unified transform method and long-range instabilities for the inhomogeneous time-dependent Schrödinger equation on the quarter-plane, *Math. Ann.* (2023)
[6] A. Chatziafratis, S. Kamvissis, A note on uniqueness for linear evolution PDEs posed on the quarter-plane, preprint (2023). arXiv:2401.08531 [math.AP]
[7] A. Chatziafratis, E.C. Aifantis, A. Carbery, A.S. Fokas, Integral representations for the double-diffusion system on the half-line, *Z. Angew. Math. Phys.* 75 (2024).
[8] A. Chatziafratis, A.S. Fokas, E.C. Aifantis, On Barenblatt's pseudoparabolic equation with forcing on the half-line via the Fokas method, Z Angew Math Mech 104 (2024).
[9] A. Chatziafratis, L. Grafakos, S. Kamvissis, Long-range instabilities for linear evolution PDEs on semi-unbounded domains via the Fokas method, *Dyn. PDE* (2024).

[10] A. S. Fokas, A unified transform method for solving linear and certain nonlinear PDEs, *Proc. Roy. Soc. London Ser. A* 453, 1411-1443 (1997).
[11] A. S. Fokas, On the integrability of linear and nonlinear PDEs, *J. Math. Phys.* 41, 4188-4237 (2000).
[12] A. S. Fokas, A new transform method for evolution partial differential equations, *IMA J. Appl. Math.* 67(6), 559-590 (2002).
[13] A. S. Fokas, A Unified Approach to Boundary Value Problems, CBMS-NSF Series Appl Math 78, SIAM (2008).
[14] A. S. Fokas, Lax Pairs: A novel type of separability (invited paper), *Inverse Problems* 25, 1-44 (2009).
[15] A. S. Fokas, E. A. Spence, Synthesis, as opposed to separation, of variables, *SIAM Review* 54 (2012).
[16] B. Deconinck, T. Trogdon, V. Vasan, The method of Fokas for solving linear partial differential equations, *SIAM Review* 56(1), 159-186 (2014).
[17] A. S. Fokas, B. Pelloni (editors), Unified Transform for Boundary Value Problems: Applications and Advances, SIAM, Philadelphia, PA (2015).
[18] A.S. Fokas and E. Kaxiras, Modern Mathematical Methods for Computational Sciences and Engineering, World Scientific (2022).